\documentclass[12pt]{article}
\usepackage{amsmath, amssymb, latexsym}
\renewcommand{\baselinestretch}{1.15}
 \newtheorem{theorem}{Theorem}[section]
\newtheorem{lemma}[theorem]{Lemma}
\newtheorem{proposition}[theorem]{Proposition}
\newtheorem{corollary}[theorem]{Corollary}
\newtheorem{remark}[theorem]{Remark}

\newtheorem{definition}[theorem]{Definition}

\newenvironment{proof}[1][Proof]{\textbf{#1.} }{\ \rule{0.5em}{0.5em}}

\renewcommand{\title}[1]
{\thispagestyle{empty}
\begin{center}
{\Large \bf #1}
\end{center}}

\newcommand{\authors}[1]
{\begin{center}
\renewcommand{\thefootnote}{\fnsymbol{footnote}}
\setcounter{footnote}{3} {\sc #1 }
\end{center}}

\newcommand{\ack}[1]{\footnote{#1}}

\newcommand{\address}[1]
{\vskip 5ex
\renewcommand{\baselinestretch}{1}
\footnotesize \normalsize
 #1 \\
}

\begin{document}

\title{Characterization of rearrangement invariant spaces with fixed points for the Hardy--Littlewood maximal operator}

\authors{
Joaquim Mart\'{\i}n\ack{Research partially supported by Grants BFM2001-3395,   
2001SGR00069 and by Programa Ram\'on y Cajal (MCYT).}
and
Javier Soria\ack{Research partially supported by Grants BFM2001-3395,   
2001SGR00069.} }

\bigskip

 {\narrower\noindent \textbf{Abstract.} \small{We characterize the rearrangement invariant spaces for which there exists a non-constant fixed point, for the Hardy--Littlewood maximal operator (the case  for the spaces $L^p(\mathbb{R}^{n})$ was first considered in \cite{Ko}). The main result that we prove is that the space $L^{\frac{n}{n-2},\infty }(\mathbb{R}^{n})\cap
L^{\infty }(\mathbb{R}^{n})$ is  minimal among those having this property.}\par}

\bigskip

\section{Introduction}

The centered Hardy--Littlewood maximal operator $\mathcal{M}$ is defined on
the Lebesgue space $L_{\rm loc}^{1}(\mathbb{R}^{n})$ by
\begin{equation*}
\mathcal{M}f(x)=\sup_{r>0}\frac{1}{\left| B_{r}\right| }\int_{B_{r}}\left|
f(x-y)\right| dy,
\end{equation*}
where $\left| B_{r}\right| $ denotes the measure of the Euclidean ball $B_{r}
$ centered at the origin of $\mathbb{R}^{n}.$

In this paper we study the existence of non--constant fixed points of the
maximal operator $\mathcal{M}$ (i.e., $\mathcal{M}f=f$) in the framework of
the rearrangement invariant (r.i.) functions spaces (see Section 2 below).  We will use some of the estimates proved in \cite{Ko}, where the case $L^p(\mathbb{R}^{n})$ was studied, and show that they can be sharpened to obtain all the rearrangement invariant norms with this property (in particular we extend Korry's result to the end point case $p=n/(n-2)$,  where the weak--type spaces have to considered.) The main argument behind this problem is the existence of a minimal space $L^{\frac{n}{n-2},\infty }(\mathbb{R}^{n})\cap
L^{\infty }(\mathbb{R}^{n})$ contained in all the r.i. spaces with the fixed point property.

\section{Background on Rearrangement Invariant Spaces}

Since we work in the context of rearrangement invariant spaces it will be
convenient to start by reviewing some basic definitions about these spaces.

A rearrangement invariant space $X=X(\mathbb{R}^{n})$ (r.i. space) is a
Banach function space on $\mathbb{R}^{n}$ endowed with a norm $\left\| \cdot
\right\| _{X(\mathbb{R}^{n})}$ such that
\begin{equation*}
\left\| f\right\| _{X(\mathbb{R}^{n})}=\left\| g\right\| _{X(\mathbb{R}^{n})}
\end{equation*}
whenever $f^{\ast }=g^{\ast }.$ Here $f^{\ast }$ stands for the
non-increasing rearrangement of $f$, i.e., the non-increasing,
right-continuous function on $\left[ 0,\infty \right) $ equimeasurable with $
f.$

An r.i. space $X(\mathbb{R}^{n})$ has a representation as a function space on
$\bar{X}(0,\infty )$ such that
\begin{equation*}
\left\| f\right\| _{X(\mathbb{R}^{n})}=\left\| f^{\ast }\right\| _{\bar{X}
(0,\infty )}.
\end{equation*}

Any r.i. space is characterized by its \textbf{fundamental function}
\begin{equation*}
\phi _{X}(s)=\left\| \chi _{E}\right\| _{X(\mathbb{R}^{n})}
\end{equation*}
(here $E$ is any subset of $\mathbb{R}^{n}$ with $\left| E\right| =s$) and
the \textbf{fundamental indices}
\begin{equation*}
\overline{\beta }_{X}=\inf\limits_{s>1}\dfrac{\log M_{X}(s)}{\log s}\text{ \
\ and \ \ }\underline{\beta }_{X}=\sup\limits_{s<1}\dfrac{\log M_{X}(s)}{
\log s},
\end{equation*}
where
\begin{equation*}
M_{X}(s)=\sup\limits_{t>0}\dfrac{\phi _{X}(ts)}{\phi _{X}(t)},\text{ }s>0.
\end{equation*}
It is well known that

\begin{equation*}
0\leq \underline{\beta }_{X}\leq \overline{\beta }_{X}\leq 1.
\end{equation*}
(We refer the reader to \cite{BS} for further information about r.i. spaces.)

\section{Main result}

Before formulating our main result, it will be convenient to start with the
following remarks (see \cite{Ko}):

\begin{remark} {\rm 
By the Lebesgue's differentiation theorem one easily obtains that
\begin{equation*}
\left| f(x)\right| \leq \mathcal{M}f(x)\text{ a.e. }x\in \mathbb{R}^{n},
\end{equation*}
thus $f$ is a fixed point of $\mathcal{M}$, if and only if $f$ is positive
and
\begin{equation*}
\frac{1}{\left| B(x,r)\right| }\int_{B(x,r)}f(y)dy\leq f(x)\text{ a.e. }x\in
\mathbb{R}^{n},
\end{equation*}
or equivalently $f$ is a positive \textbf{super--harmonic} function (i.e. $
\triangle f\leq 0,$ where $\triangle $ is the Laplacian operator).}
\end{remark}

\begin{remark}{\rm
If $f$ is a non-constant fixed point of $\mathcal{M},$ and $\varphi\ge 0 $
belongs to the Schawrtz class $\mathcal{S}(\mathbb{R}^{n}),$ with $\int_{
\mathbb{R}^{n}}\varphi (x)dx=1,$ then \ the function $f_{t}(x)=(f\ast \varphi
_{t})(x)$, with $\varphi _{t}(x)=t^{-n}\varphi (x/t)$ is also a non-constant
fixed point of $\mathcal{M}$ which belongs to $\mathcal{C}^{\infty }(\mathbb{
R}^{n})$ (notice that using the Lebesgue differentiation theorem, there
exists some $t>0$ such that $f_{t}$ is non-constant, since $f$ is non-constant). In particular if $X(\mathbb{R}^{n})$ is an r.i. space and $f\in X(
\mathbb{R}^{n})$ is a non-constant fixed point of $\mathcal{M},$ since $
\mathcal{S}(\mathbb{R}^{n})\subset L^{1}(\mathbb{R}^{n})\cap L^{\infty }(
\mathbb{R}^{n})$ we get that $f_{t}\in X(\mathbb{R}^{n})\cap \mathcal{C}
^{\infty }(\mathbb{R}^{n})$ is a non-constant fixed point of $\mathcal{M}.$}
\end{remark}

\begin{remark} {\rm 
Using the theory of weighted inequalities for $\mathcal{M}$ (see  \cite{GcRf}), if $\mathcal{M} f=f$, in particular $f\in A_1$ (the Muckenhoupt weight class), and hence $f(x)\,dx$ defines a doubling measure. Hence, $f\notin L^1(\mathbb{R}^{n})$. Also, using the previous remark we see that if $f\in L^p(\mathbb{R}^{n})$ is a fixed point, then $f\in L^q(\mathbb{R}^{n})$, for all $p\le q\le\infty$. }
\end{remark}

\begin{definition}
Given an r.i. space $X(\mathbb{R}^{n})$, we define
\begin{equation*}
D_{I_{2}}(X(\mathbb{R}^{n}))=\left\{ f\in L^{0}(\mathbb{R}^{n}):\left\|
I_{2}f\right\| _{X(\mathbb{R}^{n})}<\infty \right\},
\end{equation*}
where $I_{2}$ is the Riesz potential,
\begin{equation*}
(I_{2}f)(x)=\int_{\mathbb{R}^{n}}\left| x-y\right| ^{2-n}f(y)dy.
\end{equation*}
\end{definition}

It is not hard to see that the space $D_{I_{2}}(X(\mathbb{R}^{n}))$ is either trivial or is  the
largest r.i. space which is mapped by $I_{2}$ into $X(\mathbb{R}^{n}),$ and
is also related with the theory of the optimal Sobolev embeddings (see \cite
{petal} and the references quoted therein).

\begin{theorem}
Let $X(\mathbb{R}^{n})$ be an r.i. space. The following statements are
equivalent:

\begin{enumerate}
\item  There is a non-constant fixed point $f\in X(\mathbb{R}^{n})$ of $
\mathcal{M}$

\item  $n\geq 3$ and $\left| x\right| ^{2-n}\chi _{\{x:\left| x\right|
>1\}}(x)\in X(\mathbb{R}^{n}).$

\item  $n\geq 3$ and $\chi _{\lbrack 0,1]}(t)+t^{2/n-1}\chi _{\lbrack
1,\infty )}(t)\in \bar{X}(0,\infty ).$

\item  $n\geq 3$ and $(L^{\frac{n}{n-2},\infty }(\mathbb{R}^{n})\cap
L^{\infty }(\mathbb{R}^{n}))\subset X(\mathbb{R}^{n}).$

\item  $n\geq 3$ and $D_{I_{2}}(X(\mathbb{R}^{n}))\neq \left\{ 0\right\} .$
\end{enumerate}
\end{theorem}

\begin{proof}
$(1\rightarrow 2)$ Since if $n=1$ or $n=2,$ the only positive super--harmonic
functions are the constant functions (see \cite[Remark 1, p.\ 210]{LiLo}),
necessarily $n\geq 3.$ Moreover, it is proved in \cite{Ko}  that, if $f$ $\in
$ $\mathcal{C}^{\infty }(\mathbb{R}^{n})$ is a non-constant fixed point of $
\mathcal{M},$ then
\begin{equation*}
f(x)\geq c\left| x\right| ^{2-n}\chi _{\{x:\left| x\right| >1\}}(x).
\end{equation*}
Since   $f\in X(\mathbb{R}^{n}),$  then $\left|
x\right| ^{2-n}\chi _{\{x:\left| x\right| >1\}}(x)\in X(\mathbb{R}^{n}).$

\medskip
$(2\rightarrow 3)$ Since if $\left| x\right| ^{2-n}\chi _{\{x:\left| x\right|
>1\}}(x)\in X(\mathbb{R}^{n})$, then
\begin{equation*}
F(x)=\chi _{\{x:\left| x\right| \leq 1\}}(x)+\left| x\right| ^{2-n}\chi
_{\{x:\left| x\right| >1\}}(x)\in X(\mathbb{R}^{n}).
\end{equation*}
An easy computation shows that
\begin{equation*}
F^{\ast }(t)\simeq \chi _{\lbrack 0,1]}(t)+t^{2/n-1}\chi _{\lbrack 1,\infty
)}(t).
\end{equation*}

\medskip
$(3\rightarrow 4)$ Since $f\in (L^{\frac{n}{n-2},\infty }(\mathbb{R}^{n})\cap
L^{\infty }(\mathbb{R}^{n}))$ if and only if
\begin{equation*}
\sup_{t>0}f^{\ast }(t)W(t)<\infty,
\end{equation*}
where $W(t)=\max (1,t^{1-2/n}),$ we have that
\begin{equation*}
f^{\ast }(t)\leq \left\| f\right\| _{L^{\frac{n}{n-2},\infty }(\mathbb{R}
^{n})\cap L^{\infty }(\mathbb{R}^{n})}W^{-1}(t)
\end{equation*}
and since $W^{-1}(t)=\chi _{\lbrack 0,1]}(t)+t^{2/n-1}\chi _{\lbrack
1,\infty )}\in \bar{X}(0,\infty )$ we have that
\begin{equation*}
\left\| f\right\| _{X(\mathbb{R}^{n})}=\left\| f^{\ast }\right\| _{\bar{X}
(0,\infty )}\leq c\left\| f\right\| _{L^{\frac{n}{n-2},\infty }(\mathbb{R}
^{n})\cap L^{\infty }(\mathbb{R}^{n})}
\end{equation*}
with $c=\left\| W^{-1}\right\| _{\bar{X}(0,\infty )}.$

\medskip
$(4\rightarrow 5)$ Since (see \cite{Sw} and \cite{BR})
\begin{equation*}
\left( I_{2}f\right) ^{\ast }(t)\leq c_{1}\left(
t^{2/n-1}\int_{0}^{t}f^{\ast }(s)ds+\int_{t}^{\infty }f^{\ast
}(s)s^{2/n-1}ds\right) \leq c_{2}\left( I_{2}f^{0}\right) ^{\ast }(t)
\end{equation*}
where $f^{0}(x)=f^{\ast }(c_{n}\left| x\right| ^{n})$, $c_{n}=$ measure of
the unit ball in $\mathbb{R}^{n}.$ (Observe that $\left( f^{0}\right) ^{\ast
}=f^{\ast }).$ Rewriting the middle term in the above inequalities, using
Fubini's theorem, we get
\begin{equation*}
\left( I_{2}f\right) ^{\ast }(t)\leq d_{1}\left( \frac{n}{n-2}
\int_{t}^{\infty }f^{\ast \ast }(s)s^{2/n-1}ds\right) \leq d_{2}\left(
I_{2}f^{0}\right) ^{\ast }(t),
\end{equation*}
where $f^{**}(t)=t^{-1}\int_0^t f^*(s)\,ds$. 
Thus, $f\in D_{I_{2}}(X(\mathbb{R}^{n}))$ if and only if
\begin{equation}
\left\| \int_{t}^{\infty }f^{\ast \ast }(s)s^{2/n-1}ds\right\| _{\bar{X}
(0,\infty )}<\infty .  \label{f:uno}
\end{equation}
Since
\begin{equation*}
F(t)=\int_{t}^{\infty }\chi _{\lbrack 0,1]}^{\ast \ast
}(s)s^{2/n-1}ds=c(\chi _{\lbrack 0,1]}(t)+t^{2/n-1}\chi _{\lbrack
1,\infty )}(t))
\end{equation*}
is a decreasing function, and
\begin{equation*}
F^{0}(x)=F(c_{n}\left| x\right| ^{n})\simeq \left( \chi _{\{x:\left|
x\right| \leq 1\}}(x)+\left| x\right| ^{2-n}\chi _{\{x:\left| x\right|
>1\}}(x)\right) \in L^{\frac{n}{n-2},\infty }(\mathbb{R}^{n})\cap L^{\infty }(
\mathbb{R}^{n})
\end{equation*}
we get that $\chi _{\lbrack 0,1]}^{\circ }\in D_{I_{2}}(X(\mathbb{R}^{n}))$.

Another argument to prove this part is the following:

Since, if $n\geq 3$ (see \cite[Theorem 4.18, p.\ 228]{BS})
\begin{equation*}
I_{2}:L^{1}(\mathbb{R}^{n})\rightarrow L^{\frac{n}{n-2},\infty }(\mathbb{R}
^{n})\text{ and }I_{2}:L^{\frac{n}{2},1}(\mathbb{R}^{n})\rightarrow
L^{\infty }(\mathbb{R}^{n})
\end{equation*}
is bounded, we have that
\begin{equation*}
I_{2}:(L^{1}(\mathbb{R}^{n})\cap L^{\frac{n}{2},1}(\mathbb{R}^{n}))\rightarrow
(L^{\frac{n}{n-2},\infty }(\mathbb{R}^{n})\cap L^{\infty }(\mathbb{R}
^{n}))\subset X(\mathbb{R}^{n})
\end{equation*}
is bounded, and hence $L^{1}(\mathbb{R}^{n})\cap L^{\frac{n}{2},1}(\mathbb{R}
^{n})\subset D_{I_{2}}(X(\mathbb{R}^{n}))$.

\medskip
$(5\rightarrow 1)$ Since $n\geq 3,$ we can use the classical formula of
potential theory (see \cite[p.\ 126]{St})
\begin{equation*}
-h=\triangle (I_{2}h)
\end{equation*}
to conclude that there is a positive function $f=I_{2}\chi _{\lbrack
0,1]}^{\circ }\in X(\mathbb{R}^{n}).$ Then $0\leq f_{t}=I_{2}(\chi _{\lbrack
0,1]}^{\circ }\ast \varphi _{t})\in X(\mathbb{R}^{n})\cap \mathcal{C}
^{\infty }(\mathbb{R}^{n})$ and $\triangle f_{t}\leq 0.$
\end{proof}

\bigbreak 

We now consider particular  examples, like the Lorentz spaces:

\begin{corollary}
Let $1\leq p< \infty ,$ and assume $\Lambda^p (\mathbb{R}^{n},w)$ is a Banach space (i.e., $w\in B_p$ if $1<p<\infty$ or $p\in B_{1,\infty}$ if $p=1$, see \cite{cpss}). Then, there exists a non-constant
function $f\in \Lambda^p (\mathbb{R}^{n},w)$ such that $\mathcal{M}(f)=f$ if and
only if $n\geq 3$ and 
$$
\int_1^\infty\frac{w(t)}{t^{p(1-2/n)}}\,dt<\infty .
$$
In particular, this condition always holds, for $p>1$ and  $n$ large enough.
\end{corollary}

\begin{proof}
 The integrability condition follows by using   the previous theorem. Now, if $w\in B_p$, then there exists an $\varepsilon>0$ such that $w\in B_{p-\varepsilon}$, and hence, it suffices to take $n>2/\varepsilon$. Observe that if $w=1$ and $p=1$, then $\Lambda^1(\mathbb{R}^{n},w)=L^1(\mathbb{R}^{n})$, which does not have the fixed point property for any dimension $n$.
\end{proof}

\begin{corollary}
 Let $1\leq p, q\leq \infty $ (if $p=1$ we only consider $q=1$). Then, there exists a non-constant
function $f\in L^{p,q}(\mathbb{R}^{n})$ such that $\mathcal{M}(f)=f$ if and
only if $n\geq 3$ and 
$$
\begin{cases}
     n/(n-2)<p\leq \infty & \text{ } \\
       \quad\text{ or } & \\

     p=n/(n-2) \text{ and } q=\infty. &
\end{cases}
$$
\end{corollary}

\begin{corollary}
(See \cite{Ko}) Let $1\leq p\leq \infty .$ There exists a non-constant
function $f\in L^{p}(\mathbb{R}^{n})$ such that $\mathcal{M}(f)=f$ if and
only if $n\geq 3$ and $n/(n-2)<p\leq \infty .$
\end{corollary}

It is interesting  to know when given an r.i. space $X(\mathbb{R}^{n})$
, the space $D_{I_{2}}(X(\mathbb{R}^{n}))$ is not trivial, or equivalently
\begin{equation}
\overline{D_{I_{2}}(X(\mathbb{R}^{n}))}:=\left\{ f\in L^{0}([0,\infty
)):\left\| \int_{t}^{\infty }f^{\ast \ast }(s)s^{2/n-1}ds\right\| _{\bar{X}
(0,\infty )}<\infty \right\}   \label{spaceY}
\end{equation}
is not trivial. This will be done in terms of the fundamental indices of $X
$. We start by computing the fundamental function of $D_{I_{2}}(X(\mathbb{R}
^{n})).$

\begin{lemma}
Let $X$ be an r.i. space on $\mathbb{R}^{n},$ $n\geq 3.$ Let $Y$ be given by (
\ref{spaceY}). Then
\begin{equation*}
\phi _{Y}(s)\simeq s^{n/2}\left\| P_{1-2/n}\chi _{\lbrack 0,s]}\right\| _{X}
\end{equation*}
where $P_{1-2/n}f(t)=t^{2/n-1}\int_{0}^{t}f(s)s^{-2/n}ds.$
\end{lemma}

\begin{proof}
\begin{eqnarray*}
s^{n/2}P_{1-2/n}\chi _{\lbrack 0,s]}(t) &\simeq &s^{n/2}(\chi _{\lbrack
0,s]}(t)+\left( \frac{s}{t}\right) ^{1-2/n}\chi _{\lbrack s,\infty )}(t)) \\
&\simeq &\int_{t}^{\infty }\chi _{\lbrack 0,s]}^{\ast \ast }(r)r^{2/n-1}dr.
\end{eqnarray*}
 \end{proof}

\begin{theorem}
Let $X$ be an r.i. space on $\mathbb{R}^{n},$ $n\geq 3.$ Let $Y$ be given by (
\ref{spaceY}). Then

\begin{enumerate}
\item  If $\overline{\beta }_{X}<1-2/n$, then $Y\neq \{0\}.$

\item  If $Y\neq \{0\}$ then $\underline{\beta }_{X}\leq 1-2/n$.
\end{enumerate}
\end{theorem}

\begin{proof}
$1.)$ Let $\chi _{r}=\chi _{\lbrack 0,r]}$. Then
\begin{equation*}
P_{1-2/n}\chi _{r}(t)=\int_{0}^{1}\chi _{r}(\xi t)\frac{d\xi }{\xi ^{n/2}}
\leq c\sum_{k=0}^{\infty }2^{-k(1-n/2)}\chi _{2^{k}r}(t).
\end{equation*}
Thus
\begin{equation*}
\left\| P_{1-2/n}\chi _{r}\right\| _{X}\leq c\sum_{k=0}^{\infty
}2^{-k(1-n/2)}\phi _{X}(2^{k}r)\leq c\phi _{X}(r)\sum_{k=0}^{\infty
}2^{-k(1-n/2)}M_{X}(2^{k}).
\end{equation*}
Let $\varepsilon >0$ be such that $\overline{\beta }_{X}+\varepsilon <1-2/n.$
Then by the definition of $\overline{\beta }_{X}$ it follows readily that
there is a constant $c>0$ such that
\begin{equation*}
M_{X}(2^{k})\leq c2^{k(\overline{\beta }_{X}+\varepsilon )},
\end{equation*}
and hence
\begin{equation*}
\sum_{k=0}^{\infty }2^{-k(1-n/2)}M_{X}(2^{k})\leq \sum_{k=0}^{\infty
}2^{-k(1-n/2-\overline{\beta }_{X}-\varepsilon )}<\infty ,
\end{equation*}
which implies that $\chi _{r}\in Y.$

$2.)$ Since $Y\neq \{0\}$ if and only if $\left\| P_{1-2/n}\chi _{\lbrack
0,1]}\right\| _{X}<\infty $ and
\begin{equation}
\sup_{t>0}\left( P_{1-2/n}\chi _{\lbrack 0,1]}\right) ^{\ast \ast }(t)\phi
_{X}(t)\leq \left\| P_{1-2/n}\chi _{\lbrack 0,1]}\right\| _{X}<\infty,
\label{f:dos}
\end{equation}
and easy computations show that (\ref{f:dos}) implies that
\begin{equation}
1\leq \sup_{t\geq 1}\frac{\phi _{X}(t)}{t^{1-2/n}}=c<\infty  , \label{f:tres}
\end{equation}
then, by (\ref{f:tres})
\begin{eqnarray*}
M_{X}(a) &=&\max \left( \sup\limits_{t\geq 1/a}\dfrac{\phi _{X}(ta)}{\phi
_{X}(t)},\sup\limits_{t<1/a}\dfrac{\phi _{X}(ta)}{\phi _{X}(t)}\right)  \\
&=&\max \left( \sup\limits_{t\geq 1/a}\dfrac{\phi _{X}(ta)}{(at)^{1-2/n}}
\frac{(at)^{1-2/n}}{\phi _{X}(t)},\sup\limits_{t<1/a}\dfrac{\phi _{X}(ta)}{
\phi _{X}(t)}\right)  \\
&\simeq &\max \left( a^{1-2/n}\sup\limits_{t\geq 1/a}\frac{t^{1-2/n}}{\phi
_{X}(t)},\sup\limits_{t<1/a}\dfrac{\phi _{X}(ta)}{\phi _{X}(t)}\right) .
\end{eqnarray*}
Thus, if $a<1,$ using again (\ref{f:tres}) we get
\begin{equation*}
M_{X}(a)\geq a^{1-2/n}\sup\limits_{t\geq 1/a}\frac{t^{1-2/n}}{\phi _{X}(t)}
\geq a^{1-2/n}
\end{equation*}
which implies that
\begin{equation*}
\underline{\beta }_{X}\leq 1-2/n.
\end{equation*}
\end{proof}

Let us see that the converse in the previous theorem is not true.

\begin{proposition}
There are rearrangement invariant spaces $X$ such that

\begin{enumerate}
\item  $Y\neq \{0\}$ and $\overline{\beta }_{X}\geq 1-2/n.$

\item  $Y=\{0\}$ and $\underline{\beta }_{X}<1-2/n$.
\end{enumerate}
\end{proposition}

\begin{proof}
Let $\varphi (t)=t^{a}\chi _{\lbrack 0,1]}(t)+t^{b}\chi _{\lbrack 1,\infty
)}(t)$, with $0\leq a,b\leq 1.$ Let
\begin{equation*}
X=\left\{ f\in L^{0}([0,\infty )):\sup_{t>0}f^{\ast \ast }(t)\varphi
(t)<\infty \right\}.
\end{equation*}
Since $\varphi $ is a quasi-concave function, we have that
\begin{equation*}
\varphi (t)=\phi _{X}(t)
\end{equation*}
and
\begin{equation*}
\underline{\beta }_{X}=\min (a,b)\text{, }\overline{\beta }_{X}=\max (a,b).
\end{equation*}
On the other hand, the space $Y$ defined by (\ref{spaceY}) is not trivial if
and only if
\begin{equation*}
b\leq 1-2/n.
\end{equation*}

Now, to prove 1) take $b\leq 1-2/n$ and $a\geq 1-2/n.$ And to see 2) take $
b>1-2/n$ and $a\leq 1-2/n.$
\end{proof}

\begin{remark}{\rm 
If we consider
\begin{equation*}
X_{0}=\left\{ f\in L^{0}([0,\infty )):\sup_{t>0}f^{\ast \ast
}(t)t^{1-2/n}(1+\log ^{+}t)<\infty \right\}
\end{equation*}
and
\begin{equation*}
X_{1}=\left\{ f\in L^{0}([0,\infty )):\sup_{t>0}f^{\ast \ast }(t)\frac{
t^{1-2/n}}{(1+\log ^{+}t)}<\infty \right\}
\end{equation*}
then $\underline{\beta }_{X_{i}}=$ $\overline{\beta }_{X_{i}}=1-2/n,$   $
Y_{0}=\{0\}$ and $Y_{1}\neq \{0\}$.}
\end{remark}

\begin{remark}{\rm 
It was proved in \cite{Ko} that if we consider the strong maximal function (i.e, the maximal operator associated to centered intervals in $\mathbb{R}^{n}$), then there were no fixed points in any $L^p(\mathbb{R}^{n})$ space, regardless of the dimension. The same argument works to show that $L^p(\mathbb{R}^{n})$ cannot be replaced by any different r.i. space. Also, if we study this question for other  kind of sets, like, e.g., Buseman--Feller differentiation bases (see \cite{gu}), then the only possible fixed points are the constant functions. This observation applies to any non-centered maximal operator (with respect to balls, cubes, etc.)
}
\end{remark}

\bigbreak

\address{
\noindent
Joaquim Mart\'{\i}n\\ Dept. Mathematics
\\ U. Aut\`onoma de Barcelona\\ E-08193 Bellaterra (Barcelona),
 SPAIN\ \ \ {\sl E-mail:} 
 {\tt jmartin@mat.uab.es}

\medskip
\noindent
Javier Soria\\ Dept. Appl. Math. and Analysis
\\ University of Barcelona\\ E-08071 Barcelona,
 SPAIN\ \ \ {\sl E-mail:} 
 {\tt soria@mat.ub.es}
}

\end{document}